\documentstyle[12pt]{article}
\begin{document}
\title{\bf Multiple positivity and the Riemann zeta-function.}
\author{Olga M. Katkova\\
Dept. of Math., Kharkov National University, \\
Svobody sq., 4, 61077, Kharkov, Ukraine \\
e-mail: katkova@ilt.kharkov.ua} \maketitle
\begin{abstract}
In the paper it is discussed the relations of the Riemann
$\zeta-$function to classes of generating functions of multiply
positive sequences according to Schoenberg (also called P\'olya
frequency sequences).
\end{abstract}

{\it 1991 Mathematics Subject Classification 30D15, 30D10}


\section{Introduction.}

This paper was inspired by works of G.Csordas, T.S.Norfolk,
R.S.Varga (\cite{cnv}) and D.K.Dimitrov (\cite{d}). To formulate
their results we need the following definition.

{\bf Definition 1.} A real entire function $f$ is said to be in
the Laguerre-P\'olya class, written $f\in L-P$, if $f$ can be
represented in the form
$$ f(z)=C z^n e^{-\beta z^2 + \gamma z}\prod_{k=1}^\infty (1+\alpha_k
z)e^{-\alpha_k z},$$ where $C, \gamma, \alpha_k$ are real, $\beta
\geq 0$,  $n\in {\bf N \cup \{0\}}, \sum\alpha_k^2 <\infty.$

P\'olya and Schur (\cite{polsch}) termed a real entire function
$f$ as a function of type I  in the Laguerre-P\'olya class,
written $f\in L-P I$, if one of the functions $f(z)$, $f(-z)$,
$-f(z)$ or $-f(-z)$ can be represented in the form.
\begin{equation}
\label{a1}f(z)=C z^n e^{\gamma z}\prod_{k=1}^\infty (1+\alpha_kz),
\end{equation}
where $C\ge 0, n\in {\bf N \cup \{0\}}, \gamma\ge 0, \alpha_k\ge0,
\sum \alpha_k <\infty.$

These classes were introduced by Laguerre (see \cite{lag},
pp.168-178). Laguerre proved that $f\in L-P I$ iff $f$ is uniform
limit, on compact subsets of ${\bf C}$, of polynomials with only
real nonpositive (or real nonnegative) zeros. Laguerre supposed
that $f\in L-P $ iff $f$ is uniform limit, on compact subsets of
${\bf C}$, of polynomials with only real zeros. This theorem was
proved by P\'olya later (see \cite{pol}, p.54). We will write $
f(z) =\sum\limits_{k=0}^\infty a_kz^k \in L-P^+,$ if $f\in L-P I$
and $a_k\geq 0$ for all $k=0,1,2,\ldots$

In works (\cite{cnv}), (\cite{cv}), (\cite{cvar}), (\cite{cs}),
(\cite{d}), (\cite{cy}) it is discussed the different properties
of the Laguerre-P\'olya classes and the relations of the Riemann
$\zeta-$function to these classes.

We proceed to briefly review here some of the nomenclature pertaining
to the Riemann Hypothesis. To begin with, the Riemann $\xi$-function can be defined (see \cite {tit}, p. 16, cf. \cite{pol}, p.285) by
\begin{equation}
\label{a2}\xi(s)=\frac{1}{2}s(s-1)\pi^{-s/2}\Gamma (s/2)\zeta (s),
\end{equation}
where $\zeta$ is the Riemann $\zeta$-function. It is known that
$\xi$ is an entire function of order one and of maximal type (see,
for example, \cite{tit}, p.29). The function $\xi$ satisfies the
following functional equation.
$$ \xi(1-s) =\xi (s).$$
Hence $\xi(s+1/2)$ is an even real function and
\begin{equation}
\label{a3} \xi(s+1/2)=\frac{1}{2}( s^2-\frac{1}{4})\pi^{-s/2
-1/4}\Gamma (s/2+1/4)\zeta (s+1/2)= \sum_{k=0}^\infty b_k s^{2k}.
\end{equation}
The coefficients $b_k$ can be found by the formula (\cite{tit})
\begin{equation}
\label{a4} b_k = 8\frac{2^{2k}}{(2k)!}\int_0^\infty t^{2k}
\Phi(t)dt, \quad k = 0, 1, \ldots,
\end{equation}
and
\begin{equation}
\label{a5} \Phi (t) = \sum_{n=1}^\infty (2\pi^2 n^4 e^{9t}-3\pi
n^2 e^{5t})\exp (-\pi n^2 e^{4t}).
\end{equation}
Note that Titchmarsh uses the symbol $\Xi$ for the function
$\Xi(s)= \xi(is+\frac{1}{2})$ (\cite{tit}, chapter II). The famous
Riemann Hypothesis is the statement that all the zeros of $\Xi$
are real, that is $\Xi \in L-P.$  For the various properties and
characterizations of the Riemann $\zeta-$function see, for
example, (\cite{tit}), (\cite{ed}), (\cite{post}). The change of
variables $z=s^2$ in (\ref{a3}) yields
\begin{equation}
\label{a6} \xi_1(z) = \xi (\sqrt{z}+1/2) = \sum\limits_{k=0}^\infty b_k
z^{k}.
\end{equation}
Thus, $\xi_1$ is an entire function of order $\frac{1}{2}$ and the
Riemann Hypothesis is equivalent to the statement that $\xi_1$ has
only real negative zeros, so that $\xi_1 \in L-P^+$ . The idea
sketched above belongs to P\'olya.

In works (\cite{cnv}), (\cite{cv}), (\cite{cvar}), (\cite{cs}),
(\cite{d}), (\cite{cy}) it is considered the following problem.
Does the function $\Xi$ belong to $L-P?$ In the same works authors
investigate the equivalent problem. Does $\xi_1$ belong to
$L-P^+?$ In particular, they study the following necessary
condition for $\xi_1 \in L-P^+$ (cf.(\cite{cc}), (\cite{cv}),
(\cite{polsch})):

\begin{equation}
\label{r2}
T_1(k)=b_k^2-\frac{k+1}{k}b_{k+1}b_{k-1}\geq  0,\quad k\in{\bf N},
\end{equation}
where $b_k$ are Maclaurin coefficients of $\xi_1.$ (In today's
terminology, the inequalities of (\ref{r2}) are called Turan
inequalities). In 1927 P\'olya raised the question of whether or
not the Turan inequalities (\ref{r2}) for the function $\xi_1$ are
all valid. In 1995 G.Csordas, T.S.Norfolk and R.S.Varga
(\cite{cnv}) proved inequalities (\ref{r2}) for the function
$\xi_1.$ For various extensions and inequalities related to
(\ref{r2}), we refer to (\cite{cc}) and (\cite{d}). In (\cite{cc})
T.Craven and G.Csordas investigated certain polynomial invariants
and used them to prove that a necessary condition that some
function $f \in L-P^+$ is that its coefficients satisfy the
following double Turan inequalities $$
T_2(k)=T_1(k)^2-T_1(k+1)T_1(k-1)\geq  0,\quad k\in{\bf N},$$ where
$T_1(k)$ is defined by (\ref{r2}). But in (\cite{cs}) G.Csordas
points out that the question whether or not the higher iterated
Turan inequalities $$
T_n(k)=(T_{n-1}(k))^2-T_{n-1}(k+1)T_{n-1}(k-1)\geq  0,\quad (k\geq
n\geq 2) $$ hold for functions in the Laguerre-P\'olya class
remains open. In the same work the author pays attention to one
more open problem: whether or not the double Turan inequalities
are valid for the function $\xi_1.$

In this paper we also study the problem of whether or not the
function $\xi_1 \in L-P^+.$ But our method is slightly different
from the one described above. It is based on a characteristic
property of $L-P^+,$ which was obtained by Aissen, Schoenberg,
Whitney and Edrei in~\cite{aissen}. We need some definitions and
notations.

{\bf Definition 2.} Let m be any positive integer. The sequence
$\{a_k\}_{k=0}^\infty $ is called {\it $m$-times positive (totally
positive) } if all minors of order $ \leq m$ (of any order) of the
infinite matrix

\begin {equation}
 \left\|
  \begin{array}{ccccc}
   a_0 & a_1 & a_2 & a_3 &\ldots \\
   0   & a_0 & a_1 & a_2 &\ldots \\
   0   &  0  & a_0 & a_1 &\ldots \\
   0   &  0  &  0  & a_0 &\ldots \\
   \vdots&\vdots&\vdots&\vdots&\ddots
  \end{array}
 \right\|
\label{g1}
\end {equation}
are non-negative. The multiply positive sequences (also called
P\'olya frequency sequences) were  introduced by Fekete in 1912
(see~\cite{fek}) in connection with the problem on the exact
calculation of the number of positive zeros of the real
polynomial. We will denote by $PF_m (PF_\infty) $ the class of all
$m$-times positive (totally positive) sequences. We will denote by
$SPF_m$ the class of all sequences from $PF_m$ such that all
minors of order $\leq m$ of matrix (\ref{g1}) without vanishing
rows (columns) are positive.

We need also such notion.

{\bf Definition 3.} Let $m$ be any positive integer. The sequence
$\{a_k\}_{k=0}^\infty $ is called {\it asymptotically $m$-times
positive} if there exists a positive integer $N$ such that all
minors of matrix $$ A_N := (a_{N+j-l}),\quad\ l=0,1, \ldots , m-1,
\quad j=0,1, \ldots (a_k=0 \quad \mbox{for}\quad  k<0)$$ are
nonnegative.

We will denote by $APF_m $ the class of all asymptotically
$m$-times positive sequences. Obviously, $SPF_m \subset PF_m
\subset APF_m.$

The classes of corresponding generating functions
$$
f(z)=\sum\limits_{k=0}^\infty a_kz^k
$$
are also denoted by $PF_m, SPF_m $ and $APF_m$.

The class $PF_\infty$ was completely described by Aissen,
Schoenberg, Whitney and Edrei in~\cite{aissen} (see also
\cite{tp}, p. 412):

{\bf Theorem ASWE. } {\it A function $f \in PF_\infty$ iff
$$
f(z)=C z^n e^{\gamma z}\prod_{k=1}^\infty \frac {
1+\alpha_kz}{1-\beta_kz},
$$
where $C\ge 0, n\in {\bf Z},
\gamma\ge0,\alpha_k\ge0,\beta_k\ge0,\sum(\alpha_k+\beta_k)
<\infty.$}

In connection with this result Schoenberg \cite{schoenb} formulated the
problem of characterizing functions $f\in PF_m, m<\infty,$ and investigated
zero sets of polynomials $f\in PF_m, m<\infty.$ He proved the following
two theorems.

{\bf Theorem A} (\cite{schoenb}).  { \it Let $f$ be a polynomial
of degree $n$ and let $f\in PF_m,  m<\infty.$ Then $f(z) \neq 0$
for $z \in \{z: |\arg z| < \frac{\pi m}{n+m-1}  \}.$}

{\bf Theorem B} (\cite{sch}).  { \it Let $f$ be a polynomial,
$f(0)
> 0.$ Then $$f(z) \neq 0 \quad \mbox{for} \quad z \in \{z: |\arg
z| < \frac{\pi m}{m+1}  \} \Longrightarrow  f\in PF_m.$$ }

Both estimates in Theorems A and B are sharp.

The following fact is a simple corollary of Theorem ASWE and (\ref {a1}).

{\bf Theorem C.} {\it Let $f$ be an entire function. Then $f\in
PF_\infty$ iff $f\in L-P^+.$}

So, the Riemann Hypothesis is equivalent to statement that $\xi_1\in
PF_\infty.$ Since $PF_\infty=\cap_{m=1}^{\infty}PF_m$, this fact means
that $\xi_1\in PF_m$ for every $m\in {\bf N}.$ In this paper we investigate
the relations of the function $\xi_1$ to classes $APF_m$ and $PF_m$.
With the help of Theorem B it is easy to prove the following fact.

{\bf Theorem 1.} {\it $\xi_1 \in PF_{44}.$ }

{\it Proof of Theorem 1.} It is well-known (see, for example
\cite{tit}, chapter XV) that
\begin{equation}
\label{r1} \xi(s)\neq 0 \quad \mbox{for} \quad s\in \{s:
\frac{1}{2} < Re\  s \leq 1, \quad 0 \leq Im\  s \leq 14 \}.
\end{equation}
By (\ref{a6}) this implies that the function $\xi_1$ has no zeros
in the angle $\{ z: |\arg z| < \frac{43 \pi}{44} \}.$ It follows
from (\ref{a4}), (\ref{a5}) and (\ref{a6}) that $\xi_1(0) > 0.$
Since $\xi_1$ is an entire function of order $\frac{1}{2}$ by the
Hadamard theorem (see, for example, \cite{lev}, p. 24) we have
$$\xi_1(z) = C \prod_{k=1}^\infty (1-\frac{z}{z_k}),  $$ where
$C>0,\quad \sum\limits_{k=1}^\infty \frac{1}{|z_k|} < \infty$ and
$|\pi - \arg z_k | \leq \frac{\pi}{44}.$ Consider the polynomials
$P_N(z)= C \prod\limits_{k=1}^N (1-\frac{z}{z_k}).$ For every
$N\in {\bf N}$ the polynomial $P_N$ has no zeros in the angle $\{
z: |\arg z| < \frac{43 \pi}{44} \}.$ So by Theorem B we have $P_N
\in PF_{44},\   N\in {\bf N}.$ The function $\xi_1$ is the uniform
limit, on compact subsets of {\bf C}, of polynomials $P_N.$ Hence
$\xi_1\in PF_{44}.$ This completes the proof of Theorem 1.

Numerous calculations (see, for example, \cite{odl}, \cite{xg})
show that the height of rectangle in (\ref{r1}) can be increased
essentially. So, the constant {\it44} in Theorem 1 can be also
increased essentially. But it is impossible to prove the Riemann
hypothesis with the help of this method, because the above
arguments use the fact that the Riemann $\zeta-$function does not
vanish in the critical strip $\{s:\frac{1}{2}<Re \  s<1\}.$

We make some comments. It follows from $\xi_1\ \in PF_2$ that

$$ b_k^2-b_{k+1}b_{k-1}\geq  0,\quad k\in{\bf N}.$$

This inequalities are similar to Turan inequalities (\ref{r2}),
but they are weaker. The double Turan inequalities mean also
nonnegativity of some determinants. More precisely,

$$T_2(k)=b_k k!
 \left\|
  \begin{array}{ccc}
   b_k k! & b_{k+1}(k+1)! & b_{k+2}(k+2)!  \\
   b_{k-1} (k-1)! & b_k k! & b_{k+1}(k+1)!  \\
   b_{k-2} (k-2)! & b_{k-1}(k-1)! & b_k k!  \\

  \end{array}
 \right\|.$$

The main results of this paper are the following theorems.

{\bf Theorem 2.} {\it $\xi_1 \in APF_m$ for all $m \in {\bf N }.$
}

Note, that the Riemann Hypothesis is equivalent to the following
statement: $\xi_1 \in PF_m$ for all $m \in {\bf N }.$

{\bf Theorem 3.} {\it For every $ m\in {\bf N}$ there exists $ n_0
\in{\bf N}$ such that for all $ n \geq n_0$ the following
inclusions hold: \\
(i) $e^{nz}\xi_1(z) \in PF_m;$ \\
(ii) $\cosh(n\sqrt{z})\xi_1(z) \in PF_m.$}

Note, that multiplication by the function $e^{nz}$ does not change
zero set of $\xi_1(z)$ and multiplication by the function $\cosh
(n\sqrt{z})$ adds only negative zeros (but growth of the product
$\cosh(n\sqrt{z})\xi_1(z)$ is the same as the growth of $\xi_1$).
Hence the Riemann Hypothesis is equivalent to each of the
statements:
$$ \exists n\in{\bf N} \   \forall m \in {\bf N} : e^{nz}\xi_1(z) \in PF_m $$
or
$$\exists n\in{\bf N}\  \forall m \in {\bf N} : \cosh(n\sqrt{z})\xi_1(z) \in
PF_m.$$

The methods of proofs of Theorems 2, 3 are similar to those used
in (\cite{goost}), (\cite{ost}), (\cite{katos}), (\cite{kat}).

\section{The reduction of the proofs of Theorems 2 and 3 to the statements on
the asymptotic behavior of some multiple integrals.}

Let $f(z) = \sum _{k=0}^\infty a_k(f) z^k, \quad a_k(f) \geq 0,$
be an entire function and let
\begin{equation}
\label{a7}A_k ^ \nu(f) :=\det \|a_{k+j-l}(f)\|_{l=0,\ldots,\nu-1;
j=0, \ldots , \nu -1} \quad \ (a_k(f)=0,\ for\ k<0).
\end{equation}
Further these determinants will play a principal role. We will
denote by $A_N(f)$ the matrix
\begin{eqnarray}
\label{a8} & A_N(f) := (a_{N+j-l}(f)), \quad l=0,1, \ldots , m-1,
\quad
j=0,1, \ldots, \quad N\in{\bf N} \quad \\
\nonumber & (a_k(f)=0 \quad \mbox{for}\quad k<0).
\end{eqnarray}
In (\cite{katos}) the following fact was proved \\
{\bf Lemma 1.} $\forall k = 0, 1, 2, \ldots \forall \nu = 1, 2,
\ldots \forall r>0$
\begin{eqnarray}
\label{f12}
 &A_k^\nu(f) \nu ! (2 \pi)^\nu r^{k \nu }= \int _{-\pi}^\pi \ldots
\int _{-\pi}^\pi  \prod _{j=1} ^\nu \left ( e^{-ik\theta_j} f
(re^{ i\theta_j }) \right ) \times \nonumber
\\
&\prod _{1\leq \alpha < \beta \leq \nu }4 \sin ^2 \frac
{\theta_\alpha -\theta _\beta }{2}d\theta_1 \ldots d\theta_\nu .
\end{eqnarray}

For the reader's convenience we present the proof of Lemma 1.

{\bf Proof. } Since $$ a_k(f) = \frac {r^{-k}}{2\pi}\int
_{-\pi}^\pi
 e^{-ik\theta} f (re^{ i\theta })d\theta, \quad k \in {\bf Z}, $$
 we have
$$ A_k^\nu (f)= \frac {r^{-k \nu}}{(2\pi)^\nu}\det \| \int
_{-\pi}^\pi
f(re^{i\theta})e^{-i(k+j)\theta}e^{il\theta}\|_{l,j=0}^{\nu -1}.
$$ By virtue of the P\'olya composition formula (\cite{PoSz},
p.48, Problem 68) for any 2 sets of functions $\psi_0, \ldots ,
\psi_{\nu - 1}$ and $\varphi_0, \ldots , \varphi_ {\nu - 1}$ we
have $$ \det \| \int _a ^b \psi _\alpha (x) \varphi_\beta (x) dx
\|_{\alpha, \beta=0}^{\nu - 1}= $$ $$ \frac {1}{\nu !} \int _a^b
\ldots \int _a^b \det \| \psi _\alpha (x_\beta) \| _{\alpha, \beta
=0 }^{\nu - 1} \det \| \varphi _\alpha (x_\beta) \| _{\alpha,
\beta =0 }^{\nu - 1} dx_0 \ldots dx_{\nu - 1}.$$ Applying this
formula to
$\psi_\alpha(\theta)=f(re^{i\theta})e^{-i(k+\alpha)\theta}$ and
$\varphi_\alpha(\theta)=e^{i\alpha\theta}$ and the formula for the
Vandermonde determinant we obtain the proof of Lemma 1.

Consider the multiple integrals
\begin{eqnarray}
\label{f11} I_k^\nu (\eta,f)= \int _{-\pi}^\pi \ldots \int
_{-\pi}^\pi Re \  \left \{ \prod _{j=1} ^\nu \left (
e^{-ik\theta_j}\frac {f (e^{\eta + i\theta_j })}{f (e^\eta)}\right
) \right\} \nonumber
\\
 \prod _{1\leq \alpha
< \beta  \leq \nu }4 \sin ^2 \frac {\theta_\alpha - \theta _\beta
}{2}d\theta_1 \ldots d\theta_\nu.
\end{eqnarray}

Since $f (r) >0$ for $r \geq 0$ and the determinants $A_k^\nu (f)
$ are real such fact follows immediately from Lemma 1.

{\bf Lemma 2.} {\it $\forall k = 0, 1, 2, \ldots \quad \forall \nu
= 1, 2, \ldots \quad \forall \eta \in {\bf R}$
$$ sign I_k^\nu (\eta, f) = sign A_k^\nu (f).$$ }

Put
\begin{equation}
\label{a10} f^1 (z) = A\xi_1(z),\quad \mbox{where} \quad A>0 \quad
\mbox{is such that} \quad f^1(0) > 1.
\end{equation}

Let
\begin{equation}
\label{a12} f^2 (z) = e^{nz} f^1(z), \quad f^3 (z) =
\cosh{(n\sqrt{z})} f^1(z).
\end{equation}

Obviously Theorem 2 and Theorem 3 are equivalent to the following
statements respectively

{\bf Theorem 2'.} {\it The function $f^1\in APF_m$ for all $m \in
{\bf N }.$}

{\bf Theorem 3'.} {\it  For every $ m\in {\bf N}$ there exists $
n_0 \in{\bf N}$ such that for all $ n \geq n_0$ the following
inclusions hold: $ f^p \in PF_m \quad (p = 2, 3). $}

Now we show that Theorem 2' is a corollary of the following fact.

{\bf Proposition 1.} {\it $ \forall\nu =1, 2, \ldots , m\quad
\exists N(\nu)\in {\bf N}\ \forall k \geq N(\nu)\quad \exists \eta
= \eta(k, \nu) > 0  \Longrightarrow I_k^\nu(\eta, f^1)>0.$}

Let us deduce Theorem 2' from Proposition 1. By Proposition 1 and
Lemma 2 $$ \forall m \in {\bf N}\ \exists N(m)\in {\bf N}\ \forall
k\ge N(m)  \forall \nu = 1, 2, \ldots , m  \Longrightarrow
A_k^\nu(f^1)>0,$$ that is all minors of order $\nu=1, 2, \ldots ,
m$ composed of consecutive rows and consecutive columns of matrix
$A_N(f^1)$ (defined by (\ref{a7})) are positive. The statement on
$f^1\in APF_m$ means nonnegativity of all minors of matrix
$A_N(f^1)$. As a matter of fact these minors are positive. It
follows from such result of Schoenberg (\cite{schoenb}).

{\bf Theorem D.} {\it If all minors of order $\nu=1, 2, \ldots ,
m$ composed of consecutive rows and consecutive columns of matrix
$A$ are positive then all minors of $A$ of order $\nu=1, 2, \ldots
, m$ are positive.}

Let
\begin{equation}
\label{a13} f^2_\varepsilon (z) = e^{z} f^1(\varepsilon z), \quad
\varepsilon = \frac{1}{n},
\end{equation}
\begin{equation}
\label{a14}  f^3_\varepsilon (z) = \cosh{\sqrt{z}} f^1(\varepsilon
z), \quad \varepsilon = \frac{1}{n^2}.
\end{equation}

It is clear that $$ \forall \varepsilon > 0\ :\ (f(z)\in PF_m\quad
\Longleftrightarrow\ f(\varepsilon z)\in PF_m).$$

Thus, Theorem 3' follows from such fact.

{\bf Theorem 3''.} {\it $ \forall m \in {\bf N }\ \exists
\varepsilon_0 (m)>0\ \quad \forall \varepsilon,\
0<\varepsilon<\varepsilon_0(m) \quad\forall p=2,3 \Longrightarrow
\quad f_\varepsilon^p(z)\in PF_m.$}

We obtain Theorem 3'' as a corollary of such proposition.

{\bf Proposition 2.} {\it $ \forall\nu \in {\bf N}\quad \exists
N(\nu)\in {\bf N}\ \quad \forall \varepsilon,\
0<\varepsilon<1\quad \forall k\ge N(\nu) \quad \forall p=2,3
\Longrightarrow I_k^\nu(\eta, f_ \varepsilon^p)>0.$}

By Proposition 2 and Lemma 2 we have

\begin{equation}
\label{a15}\forall\nu \in {\bf N}\ \exists N(\nu)\in {\bf N}\
\forall \varepsilon,\ 0<\varepsilon<1,\ \forall k\ge N(\nu)
\forall p=2,3 \Longrightarrow A_k^\nu(f_\varepsilon^p)>0.
\end{equation}

Now we need the following lemma from (\cite{katos}) which is
similar to Theorem~D.

{\bf Lemma 3.}(\cite{katos}) {\it Let $\{a_k\}_{k=0}^\infty$ be a
sequence of positive numbers such that $\sum_0^\infty a_k<\infty$.
Consider the matrix with $\nu$ rows and infinitely many columns

$$ A_\nu=\|a_{j-i}\|_{i=0,\ldots,\nu-1;j=0,1,\ldots}\ (a_k=0,\
for\ k<0).$$ Suppose that for every $\nu=1,\ldots,m$ the matrix
$A_\nu$ satisfies the following condition: all minors of order
$\nu$, composed of consecutive columns, are positive. Then
$\{a_k\}\in PF_m$.}

Note that $$ f_\varepsilon^1 (z) \to e^z f^1(0),\quad
f_\varepsilon^2(z)\to \cosh {(\sqrt{z})}f^1(0)$$ uniformly on any
compact set of ${\bf C}$-plane. It is well-known that $e^z \in
SPF_\infty.$  Since the matrix (\ref{g1}) of coefficients of $
\cosh (\sqrt{z})$ is the submatrix of the matrix (\ref{g1}) of
coefficients of $e^z$ the function $ \cosh (\sqrt{z})\in
SPF_\infty.$ Therefore $$\forall \nu \in {\bf N} \quad \exists
\varepsilon (\nu) >0 \quad \forall \varepsilon,
0<\varepsilon<\varepsilon_\nu \quad \forall k =0,1, \ldots, N(\nu)
-1 \quad $$ $$\forall p= 2, 3 \Longrightarrow A_k^\nu
(f_\varepsilon^p) >0.$$
 So, by (\ref{a15})
$$\forall \nu \in {\bf N}\quad \exists \varepsilon (\nu) >0 \quad \forall
\varepsilon, 0<\varepsilon<\varepsilon_\nu \quad \forall k =0,1,2
\ldots \quad \forall p= 2,3 \Longrightarrow A_k^\nu
(f_\varepsilon^p)>0.$$

Putting $\tilde{\varepsilon} = \min \{\varepsilon (1), \ldots,
\varepsilon (m)\}$ we obtain that the sequence of coefficients of
$f_\varepsilon^p,  p =2,3$   satisfies the assumptions of Lemma 3
and thus $f_\varepsilon^p \in PF_m, p =2,3.$ We have proved that
Theorem 3'' follows from Proposition 2.

Thus the proofs of Theorems 1 and 2 are reduced to the statements
on positivity of corresponding integrals of the form (\ref{f11}).
Let $\vec{\theta} = (\theta _1, \ldots , \theta_\nu) \in {\bf
R}^\nu .$ We will use the $l_\infty$-norm. To prove Propositions 1
and 2 for $f=f^1, f^2_\varepsilon, f^3_\varepsilon $ we write the
integral (\ref{f11}) as $$ I_k^\nu (\eta,f)=\left ( {\int \ldots
\int }_{\|\vec{\theta} \| \leq \sigma}+ {\int \ldots \int }_
{\sigma < \| \vec{\theta} \| \leq \pi }\right ) \Re \left \{ \prod
_{j=1} ^\nu \left ( e^{-ik\theta_j}\frac {f (e^{\eta + i\theta_j
})}{f (e^\eta)}\right ) \right\} \times $$
 \begin{equation}
\label{f13}
    \prod _{1\leq \alpha < \beta  \leq \nu }4 \sin ^2 \frac {\theta_\alpha -
\theta _\beta }{2}d\theta_1 \ldots d\theta_\nu = J_1(f) + J_2(f),
\end{equation}
where $\sigma =\sigma (f) >0$ and $\eta = \eta (f) >0$ will be
chosen with the help of reasoning usually applied in the
saddle-point method. We will estimate $J_1(f)$ from bellow and
$J_2(f)$ from above.

\section{The estimates from bellow \\
of the integrals $J_1(f)$.}

We need some properties of the function $\xi _1$.

{\bf Lemma 4.} {\it The function $\xi _1$ satisfies the following
conditions:\\
(a) all coefficients of $\xi _1$ are positive; \\
(b) $\xi _1$ has no zeros in the angle $| \arg z | < \pi /2 $;
\\
(c) for $\theta $, $\ | \theta | \leq \pi /2, $ we have $$
\log|\xi _1 (re^{i\theta})|= \frac {1}{4} \sqrt {r} \log r \cos
\frac {\theta }{2}\   -  \frac {1}{4} \sqrt {r}\theta \sin \frac
{\theta }{2}\   -$$   $$ \frac {1}{2} \sqrt {r} (\log 2 + 1 + \log
\pi ) \cos \frac {\theta }{2} + O(\log r),\quad r\to\infty;$$ \\
(d) for $\theta $, $| \theta | < \pi , $ we have $$ \log|\xi _1
(re^{i\theta})| \leq \frac {1}{4} \sqrt {r} \log r \cos \frac
{\theta }{2} + O(\sqrt {r} \log r);$$ \\ (e) $ \xi _1^{'} (r)/ \xi
_1(r) = \frac{1}{8}\frac{\log r}{\sqrt {r}} + O(\frac {1} {\sqrt
{r}}),\quad r \to \infty.$ }\\

{\it Proof.} In virtue of (\ref{a4}), (\ref{a5}), (\ref{a6}) the
statement (a) is obvious. The proposition (b) follows immediately
from (\ref{a6}) and (\ref{r1}). To prove (c) we note that it
follows from such well-known properties of the $\zeta-$function
(see, for example \cite{tit}, chapter I)
$$\zeta(z)=\sum_{n=1}^{\infty} \frac{1}{n^z},\quad Re\  z >1$$ and
$$\zeta(z)=\prod_{p}\left(1-\frac{1}{p^z}\right)^{-1}, \quad Re \
z>1 $$(here and further we denote prime numbers by $p$) that
$$|\zeta(z)|\leq C, \quad Re\ z\geq 2$$ and $$|\zeta(z)|^{-1}\leq
C,\quad Re\ z\geq 2.$$ (Here and further we denote positive
constants by $C$ without indexes). Then from (\ref{a3}) with the
help of a well-known Stirling's Formula (see, for example
\cite{ing}, chapter III) $$ \log \Gamma (z)= z\log
z-z-\frac{1}{2}\log z+\frac{1}{2}\log{(2\pi)}+O(\frac{1} {|z|}),
\quad |z|\to \infty, \quad |\arg z|\leq\pi-\delta,$$ we obtain
(c).

Since all zeros of $\xi-$function lie in the strip $\{z:0<Re \
z<1\}$ (see, for example \cite{tit}, chapter III), than by
(\ref{a6}) the number of zeros of the function $\xi_1$ in every
angle $\{z:|\arg z|<\psi\}, \ 0<\psi\leq\pi,$ is finite. Hence,
the indicator of  the function $\xi_1$ is trigonometric. So,
 the inequality (d) is a well-known property of the indicator of an entire
function (see, for example, \cite{lev}, p.71).

Let us prove (e). It follows from the formula (see, for example,
\cite{tit}, chapter~I)
$$-\frac{\zeta^{'}(z)}{\zeta(z)}=\sum_{m=1}^{\infty}\sum_p\frac{\log
p} {p^{mz}}, \quad Re\  z >1,$$  that
$$\left|\frac{\zeta^{'}(r)}{\zeta(r)}\right|= O(1),\quad r\to
\infty.$$

By (\ref{a3}) and the formula (see, for example, \cite{ing},
chapter III) $$\frac{\Gamma^{'}(z)}{\Gamma(z)}=\log z
-\frac{1}{2z}+O(\frac{1}{|z|^2}), \quad|z|\to \infty,\quad |\arg
z|\leq\pi-\delta,$$ we have (e). Lemma 4 is proved.

Since the methods of proofs of Propositions 1, 2 are similar, we
will prove these statements at the same time. We will denote by
$b(\eta ,f)$ the function
\begin{equation}
\label{f14}
 b (\eta ,f) = \log f (e^\eta ).
\end{equation}

So,
 \begin{equation}
\label{f16}
 b(\eta ,f^1) = \log f^1 (e^\eta) = \log ( A\xi_1
 (e^{\eta }));
\end{equation}

\begin{equation}
\label{f17}
 b(\eta ,f_\varepsilon ^2) =
 e^\eta + \log f^1(\varepsilon e^\eta );
\end{equation}

\begin{equation}
\label{f18}
 b(\eta ,f_\varepsilon ^3) =  \log \cosh
 (e^{\eta/2})+ \log f^1(\varepsilon e^\eta ) .
\end{equation}

Further we need such lemma.

{\bf Lemma 5.}  {\it For $\eta \geq \eta _0$ the following
inequalities are valid}
\begin{equation}
\label{f19} b(\eta ,f^1)\leq C\eta e^{\eta/2}\leq b'(\eta ,f^1);
\end{equation}

\begin{equation}
\label{f20} b'(\eta ,f^1)\leq C\eta e^{\eta/2}\leq b(\eta ,f^1).
\end{equation}

{\it For any $\eta \geq 0$ and every $\varepsilon , 0<\varepsilon
<1$, the following inequalities hold}

\begin{equation}
\label{f21} b(\eta ,f_\varepsilon ^2)\leq C(e^\eta
+\sqrt{\varepsilon } e^{\eta/2}(\eta + \log \varepsilon ))\leq
b'(\eta , f_\varepsilon ^2);
\end{equation}

\begin{equation}
\label{f22} b'(\eta , f_\varepsilon ^2)\leq C(e^\eta
+\sqrt{\varepsilon} e^{\eta/2}(\eta + \log \varepsilon )) \leq
b(\eta , f_\varepsilon ^2);
\end{equation}

\begin{equation}
\label{f23} b(\eta ,f_\varepsilon ^3)\leq C(e^{\eta/2}
+\sqrt{\varepsilon } e^{\eta /2}(\eta + \log \varepsilon ))\leq
b^{'}(\eta , f_\varepsilon ^3);
\end{equation}

\begin{equation}
\label{f24} b^{'}(\eta , f_\varepsilon ^3)\leq C(e^{\eta /2}
+\sqrt{\varepsilon} e^{\eta/2}(\eta + \log \varepsilon )) \leq
b(\eta , f_\varepsilon ^3).
\end{equation}

{\bf Proof.} In virtue of (\ref{f16}) it is clear that
inequalities (\ref{f19}) and (\ref{f20}) are the corollaries of
the properties (c), (e) from Lemma 4.

To prove (\ref{f21})-(\ref{f24}) we note that by (c) the following
estimates are valid

\begin{equation}
\label{f25} C\sqrt {r} \log r - C \leq \log\xi _1 (r) \leq C \sqrt
{r} \log r + C,\quad r \geq 0.
\end{equation}

It follows from (e) that

\begin{equation}
\label{f26} C\sqrt {r} \log r - C \leq \frac {\xi _1^{'}(r)} {\xi
_1(r)}
 \leq C\sqrt {r} \log r + C,\quad r \geq 0.
\end{equation}

Hence, by (\ref{a10}) and (\ref{f17}) we have

$$ b(\eta ,f_\varepsilon ^2)\leq e^\eta +C\sqrt{\varepsilon }
e^{\eta/2}(\eta + \log \varepsilon )+C \leq C(e^\eta
+\sqrt{\varepsilon } e^{\eta/2}(\eta + \log \varepsilon ))$$
$$\leq C(e^\eta +C\varepsilon e^\eta(\frac {\xi_ 1^{'}(r)}{\xi_
1(r)})\left. \right|_{r=\varepsilon e^\eta}  + C)\leq C b^{'}(\eta
, f_\varepsilon ^2). $$

Inequality (\ref{f21}) is proved. The estimates (\ref{f22}),
(\ref{f23}) and (\ref{f24}) can be obtained in the same way as
(\ref{f21}). Lemma 5 is proved.

It follows from property (b) of Lemma 4, that for every
$f=f^1,f_\varepsilon ^2, f_\varepsilon ^3$ the function $f(e^{\eta
+ i \theta }) \ne 0$ in the circle $ \{ \theta \in {\bf C}:
|\theta| \leq \pi /4 \}.$ Thus, for every $f=f^1,f_\varepsilon ^2,
f_\varepsilon ^3$ and $|\theta|\le \pi /4$ we have the
decomposition

\begin{equation}
\label{b27} \log \left\{e^{-ik\theta}\frac{f (e^{\eta+i\theta})}{f
(e^\eta)}\right\}=-ik\theta+\sum_{j=1}^\infty b^{(j)}(\eta ,f)
\frac{i^j\theta^j}{j!}.
\end{equation}

By virtue of property (a) of Lemma 4 for every
$f=f^1,f_\varepsilon ^2, f_\varepsilon ^3$ the function $b'(\eta
,f)>0,$ $b'(\eta ,f^1)$ is an increasing function of $\eta ,$
$b'(\eta , f_\varepsilon ^p ), p=2,3, $ are the increasing
functions of $\eta ,$ and of $\varepsilon $. Therefore for every
$f=f^1,f_\varepsilon ^2, f_\varepsilon ^3$ there exists a positive
integer $k_0=k_0(f)$ such that for every $k\geq k_0(f)$ the
equation

\begin{equation}
\label{b28} b'(\eta , f)=k
\end{equation}
has the unique solution $\eta=\eta(k,f).$  In so doing, we will
have
$$k_0(f^1)=\max\{ [b'(0,f^1)]+1,1\};$$ $$ k_0(f_\varepsilon
^p)=\max\{\max_{0\le\varepsilon\le1}[b'(0,f_\varepsilon^p)]+1,1\}
=\max\{[b' (0,f^p_1)]+1,1\}, p=2,3.$$ Note, that $
k_0(f_\varepsilon ^p), p=2,3,$ does not depend on $\varepsilon ,
0<\varepsilon <1.$
 By the choice of $\eta$ we write (\ref{b27}) in
the form
\begin{equation}
\label{b29} \log \left\{e^{-ik\theta}\frac {f
(e^{\eta+i\theta})}{f (e^\eta)}\right\}=-\frac {1}{2}\theta^2b''
(\eta ,f)+ \tau (\theta,\eta,f),
\end{equation}
where
\begin{equation}
\label{b30} \tau (\theta,\eta,f) = \sum_{j=3}^\infty
b^{(j)}(\eta,f) \frac{i^j\theta^j}{j!}.
\end{equation}

We need the following Lemma.

{\bf Lemma 6.}{\it \  $\exists C>0\  \forall f=f^1,
f_\varepsilon^2, f_\varepsilon^3\  \exists k_0(f)\  \forall k \geq
k_0, k\in \mathbf{N},\ \forall j\in \mathbf{N}:$
\begin{equation}
\label{b31} |b^{(j)}(\eta,f)| \leq \frac {4^j j!}{\pi ^j} C k.
\end{equation}
Moreover, $ k_0(f_\varepsilon ^2)$ and $ k_0(f_\varepsilon ^3)$ do
not depend on $\varepsilon , 0<\varepsilon <1.$ }

{\bf Proof.} For every $f=f^1, f_\varepsilon ^2, f_\varepsilon ^3$
we have $f(e^{\eta + z}) \ne 0$ in the circle $ \{z\in {\bf C}:
|z| \leq \pi /4 \}.$ Applying the Schwarz formula (see, for
example, \cite{PoSz}, Problem 231) to the function $\log f
(e^{\eta +z})$ in the circle $|z| \leq \pi /4,$ we obtain $$\log f
(e^{\eta +z})=\frac{1}{2\pi}\int_{-\pi}^{\pi}\log |f (e^{\eta
+\frac{\pi}{4} e^{i\tau}})|\frac{\frac{\pi}{4}
e^{i\tau}+z}{\frac{\pi}{4} e^{i\tau}-z}d\tau +i C,$$ where $f=f^1,
f_\varepsilon ^2, f_\varepsilon ^3.$ Differentiating with respect
to $z,$ setting $z=0,$ we have $$|b^{(j)}(\eta,f)| \leq \frac {4^j
j!}{\pi ^{j+1} }\int\limits_{-\pi}^\pi |\log |f (e^{\eta +
\frac{\pi}{4}e^{i\tau }})||d\tau = $$  $$ \frac {2^{2j+1} j!}{\pi
^j }\left( \frac{1}{\pi}\int\limits_{-\pi}^\pi \log ^+ |f (e^{\eta
+ \frac{\pi}{4} e^{i\tau}})|d\tau
-\frac{1}{2\pi}\int\limits_{-\pi}^\pi \log |f (e^{\eta +
\frac{\pi}{4}e^{i\tau }})|d\tau \right) = $$ $$ \frac {2^{2j+1}
j!}{\pi ^j }\left( \frac {1}{\pi}\int\limits_{-\pi}^\pi \log ^+ |f
(e^{\eta + \frac{\pi}{4}e^{i\tau}})|d\tau - \log f (e^{\eta
})\right),$$ where $f=f^1, f_\varepsilon ^2, f_\varepsilon ^3.$
Taking into account (\ref{a10}) and (\ref{f14}), (\ref{f16}),
(\ref{f17}),  we obtain
\begin{equation}
\label{b32} |b^{(j)}(\eta,f)| \leq \frac {2^{2j+1} j!}{\pi ^{j+1}
}\int_{-\pi}^\pi \log ^+ |f (e^{\eta + \frac{\pi} {4}e^{i\tau
}})|d\tau ,
\end{equation}
where $f=f^1, f_\varepsilon ^2, f_\varepsilon ^3.$ By property (a)
of Lemma 4 for $f=f^1, f_\varepsilon ^2, f_\varepsilon ^3$ we have
\begin{equation}
\label{b33} \log ^+ |f (e^{\eta + \frac{\pi} {4}e^{i\tau }})| \leq
\log  f (e^{\eta + \frac{\pi} {4} })= b((\eta+ \frac{\pi} {4}),f).
\end{equation}

In virtue of (\ref{f19}) and (\ref{b28}) we obtain
\begin{equation}
\label{b34} b((\eta+ \frac{\pi} {4}),f^1)\leq C\eta e^{\eta/2}\leq b'(\eta ,f^1)=C k.
\end{equation}
By (\ref{f21}) and (\ref{b28}) for all $\varepsilon , 0<\varepsilon <1,$ we have
\begin{equation}
\label{b35} b((\eta+ \frac{\pi} {4}) ,f_\varepsilon ^2)\leq C(e^\eta
+\sqrt{\varepsilon } e^{\eta/2}(\eta + \log \varepsilon ))\leq
b'(\eta , f_\varepsilon ^2)= C k.
\end{equation}
where $f=f^1, f_\varepsilon ^2, f_\varepsilon ^3.$ By property (a)
of Lemma 4 for $f=f^1, f_\varepsilon ^2, f_\varepsilon ^3$ we have
\begin{equation}
\label{s33} \log ^+ |f (e^{\eta + \frac{\pi} {4}e^{i\tau }})| \leq
\log  f (e^{\eta + \frac{\pi} {4} })= b((\eta+ \frac{\pi} {4}),f).
\end{equation}

In virtue of (\ref{f19}) and (\ref{b28}) we obtain
\begin{equation}
\label{s34} b((\eta+ \frac{\pi} {4}),f^1)\leq C\eta e^{\eta/2}\leq
b'(\eta ,f^1)=C k.
\end{equation}
By (\ref{f21}) and (\ref{b28}) for all $\varepsilon , 0<\varepsilon <1,$ we have
\begin{equation}
\label{s35} b((\eta+ \frac{\pi} {4}) ,f_\varepsilon ^2)\leq
C(e^\eta +\sqrt{\varepsilon } e^{\eta/2}(\eta + \log \varepsilon
))\leq b'(\eta , f_\varepsilon ^2)= C k.
\end{equation}
By (\ref{f23}) and (\ref{b28}) for all $\varepsilon , 0<\varepsilon <1,$ we
have
\begin{equation}
\label{b36} b((\eta + \frac{\pi} {4})  ,f_\varepsilon ^3)\leq C(e^{\eta/2}
+\sqrt{\varepsilon } e^{\eta /2}(\eta + \log \varepsilon ))\leq
b^{'}(\eta , f_\varepsilon ^3)= C k.
\end{equation}

Substituting (\ref{s34})-(\ref{b36}) into (\ref{s33}) and then
into (\ref{b32}) we obtain statement of Lemma 6. Lemma 6 is
proved.

It follows from (\ref{b33}) that for $\theta, |\theta| \leq \frac
{\pi}{8}$ and for $f=f^1, f_\varepsilon ^2, f_\varepsilon ^3$ we
have $$ |\tau (\theta ,\eta , f )|\leq C k\sum _{p=3}^\infty \frac
{4^p |\theta |^p}{\pi ^p} \leq B |\theta |^3 k,$$ where $B>0.$ Let
us choose
\begin{equation}\label{b37}
\sigma = \sigma (k) =\left ( \frac {\pi}{3B\nu }\right
)^{1/3}k^{-1/3},
\end{equation}
then for $\theta , |\theta | \leq \sigma,$ we have
\begin{equation}\label{b38}
|\tau (\theta, \eta, f ) | \leq \frac {\pi}{3\nu}.
\end{equation}
Applying (\ref{b29}), (\ref{b31}) for $\vec{\theta},
\|\vec{\theta}\|\leq \sigma,$ we obtain
\begin{equation}\label{b39}
Re\  \prod _{j=1}^\nu  \left\{e^{-ik\theta _j}\frac {f
(e^{\eta+i\theta _j})}{f (e^\eta)}\right\}\geq \frac {1}{2}
e^{-\frac {\pi}{3}}\exp \left ( -C k \sum _{j=1} ^\nu \theta _j ^2
\right ) .
\end{equation}
So, with the help of (\ref{b37}) we obtain
\begin{eqnarray}
& J_1(f) \geq \left (\frac {2}{\pi} \right )^{\nu (\nu -1)} C {
\int \ldots \int }_{\|\vec{\theta} \| \leq \sigma} \exp \left ( -C
k \sum _{j=1} ^\nu \theta _j ^2 \right )  \nonumber
\\
& \prod _{1\leq \alpha < \beta \leq \nu } (\theta_ \alpha - \theta
_\beta )^2 d\theta_1 \ldots d\theta _\nu =  \nonumber \\ & \left
(\frac {2}{\pi} \right )^{\nu (\nu -1)} C k^{-\frac { \nu^2}{2} }
{ \int \ldots \int } _{\|\vec{u} \| \leq C k^{\frac {1}{6}}} \exp
\left ( -C  \sum _{j=1} ^\nu u _j ^2 \right )  \nonumber
\\
& \prod _{1\leq \alpha < \beta \leq \nu } (u_ \alpha - u _\beta
)^2 du_1 \ldots du _\nu \geq {C(\nu){k^{-\frac {\nu^2}{2}}} } ,
\label{b40}
\end{eqnarray}

where $C(\nu) >0 $ depends only on $\nu.$

\section{The estimate from above of the integrals $J_2(f)$. }

The integration domain $\{\vec{\theta}:
\sigma<\|\vec{\theta}\|\le\pi\}$ in $J_2(f)$ is contained in the
union of the domains $$
\{\vec{\theta}:|\theta_1|\le\pi,\ldots,|\theta_{j-
1}|\le\pi,\sigma<|\theta_j|\le\pi,|\theta_{j+1}|\le\pi,\ldots,|\theta_\nu|\le\pi\},
\ j=1,\ldots,\nu$$ and the integrand has a majorant $$ 2^{(\nu-
1)\nu}\prod_{j=1}^\nu\left|\frac{f (e^{\eta+i\theta_j})}{f
(e^\eta)}\right|, $$ which is symmetric with respect to
$\theta_1,\ldots,\theta_\nu.$ Therefore,
$$|J_2(f)|\le\nu2^{(\nu-
1)\nu}\int\limits_{\sigma<|\theta|\le\pi}\left| \frac{f
(e^{\eta+i\theta})}{f (e^\eta)}\right|d\theta \left(\int\limits_{-
\pi}^\pi\left|\frac{f (e^{\eta+i\theta})} {f
(e^\eta)}\right|d\theta\right)^{\nu-1}.$$

By property (a) of Lemma 4 and (\ref{a10}), (\ref{a12})we have
$$\left|\frac{f (e^{\eta+i\theta})}{f (e^\eta)}\right|<1.$$ Hence
\begin{eqnarray}
\label{b41} & |J_2(f)|\le \nu 2^{(\nu- 1)(\nu+1)}\pi^{(\nu-1)}
\left(\int\limits_{\sigma<|\theta|\le\frac{\pi }{
2}}+\int\limits_{\frac{\pi}{2}<|\theta|\le\pi}\right) \left(
\left|\frac{f (e^{\eta+i\theta})}{f (e^\eta)}\right|d\theta\right)
 \nonumber
\\ &= C(\nu)(I_1(f)+I_2(f)).
\end{eqnarray}

To estimate the integral $I_1(f)$  we note that by property (c) of
Lemma 4 for $\theta $, $\ | \theta | \leq \pi /2, $ we have $$
\log|\xi _1 (re^{i\theta})| - \log \xi _1 \leq \frac {1}{4} \sqrt
{r} \log r (\cos \frac {\theta }{2}\   - 1)- \frac {1}{4} \sqrt
{r}\theta \sin \frac {\theta }{2}\   -$$   $$ \frac {1}{2} \sqrt
{r} (\log 2 + 1 + \log \pi ) (\cos \frac {\theta }{2}-1) + C\log r
,\quad r\geq r_0.$$
 Hence for $\theta $, $\sigma\leq | \theta | \leq \pi /2, $
and $ r\geq 0$ the following inequality is valid

\begin{equation}\label{b42}
\log|\xi_1(re^{i\theta})|-\log \xi_1(r)\leq -C\sqrt {r}\log r
\sigma^2+C\log r+C.
\end{equation}

By (\ref{a10}) and (\ref{f20}) for $\theta $, $\sigma \leq |
\theta | \leq \pi /2,$ we have
\begin{eqnarray}\label{b43}
&\log|f^1(e^{\eta + i\theta})|-\log f^1(e^\eta)\le -C\eta e^{\eta
/2}\sigma^2 +C\eta \\
\nonumber & \leq -C b{'}(\eta, f ^1)\sigma^2 +C\eta,\quad \eta
\geq \eta _0.
\end{eqnarray}

Note that
\begin{equation}
\label{m5} e^{\eta +i\theta } -e^\eta \leq -2/\pi ^2 e^\eta \theta
^2,\quad |\theta|\leq \pi/2;
\end{equation}
\begin{eqnarray}
\label{m6} & \log |\cosh (e^{(\eta+i\theta)/2)})|-
\log\cosh(e^{\eta/2})\leq \\
\nonumber & e^{\frac{\eta}{2}}( \cos{\frac{\theta}{2}}-1)+ \log 2
\leq -\frac{1}{2\pi ^2}e^{\eta/2}\theta^2+\log 2, |\theta|\leq\pi.
\end{eqnarray}
Then by (\ref{a10}), (\ref{a13}), (\ref{b42}) and (\ref{f22}) for
all $\varepsilon,   0< \varepsilon <1,$ and $\eta >0$ we have
\begin{eqnarray}\label{b47}
& \log |f_\varepsilon ^2(e^{\eta+i\theta})|-\log
f_\varepsilon^2(e^\eta) )\leq \\ \nonumber &- C(e^\eta
+\sqrt{\varepsilon } e^{\eta/2}(\eta + \log \varepsilon ))\sigma^2
+ C\eta +C \leq \\ \nonumber & -C b{'}(\eta,
f_\varepsilon^2)\sigma^2 +C\eta +C, \quad \sigma\leq |\theta| \leq
\pi /2,
 \end{eqnarray}
and by (\ref{a10}), (\ref{a14}), (\ref{b42}) and (\ref{f24}) for
all $\varepsilon,  0 < \varepsilon < 1,$ and $\eta >0$ we have
\begin{eqnarray}
\label{b48} & \log |f_\varepsilon ^3(e^{\eta +i\theta})|-\log
f_\varepsilon^3(e^\eta) )\leq \\ \nonumber & -C(e^{\eta/2}
+\sqrt{\varepsilon } e^{\eta /2}(\eta + \log \varepsilon
))\sigma^2 + C\eta +C \leq \\ \nonumber &-C b{'}(\eta,
f_\varepsilon^3)\sigma^2 +C\eta +C, \quad \sigma\leq |\theta| \leq
\pi /2.
\end{eqnarray}

From (\ref{b43}), (\ref{b47}), (\ref{b48}) with the help of (\ref{b28})
and (\ref{b37}) for $\theta $, $\ \sigma <| \theta | \leq \pi /4, $ and
$f=f^1, f_\varepsilon ^2, f_\varepsilon ^3$ we obtain
\begin{equation}\label{b44}
\log|f(e^{\eta + i\theta})|-\log f(e^\eta)\le -C k^\frac
{1}{3}+C\log k+C\le -C k^\frac{1}{3}, \quad k \geq k_0(f),
\end{equation}
where $k_0(f_\varepsilon ^2), k_0( f_\varepsilon ^3)$ do not
depend on $\varepsilon$.

Hence
\begin{equation}\label{m1}
I_1 \leq C\exp ( - C(\nu)  k^\frac {1}{3}).
\end{equation}

 Now we will estimate $I_2$. By the conditions (c) and (d) of Lemma 4
 for $\theta, \frac{\pi}{2} <|\theta|< \pi$ we have
\begin{eqnarray}
\label{m2} & \log|\xi_1(re^{i\theta})|-\log \xi_1(r)\leq  \\
\nonumber & -C\sqrt {r}\log r \theta^2+o(\sqrt {r}\log r)\leq
-C\sqrt {r}\log r + C,\quad r \geq 0.
\end{eqnarray}

By (\ref{a10}) and (\ref{f20}) for $\theta $, $\frac{\pi}{2}
<|\theta | <\pi, $ we have
\begin{equation}\label{m8}
\log|f^1(e^{\eta + i\theta})|-\log f^1(e^\eta)\leq -C\eta e^{\eta
/2}\leq -C b{'}(\eta, f ^1) +C ,\quad \eta \geq \eta _0.
\end{equation}
Then by (\ref{a10}), (\ref{a13}), (\ref{m2}) and (\ref{f22}) for
all $\varepsilon, 0< \varepsilon <1,$ and $\eta >0$  we have
\begin{eqnarray}\label{m9}
&\log |f_\varepsilon ^2(e^{\eta+i\theta})|-\log
f_\varepsilon^2(e^\eta) )\leq \\ \nonumber &- C(e^\eta
+\sqrt{\varepsilon } e^{\eta/2}(\eta + \log \varepsilon ))+ C \leq
-C b{'}(\eta, f_\varepsilon^2) +C ,\quad \frac{\pi}{2} <|\theta|
<\pi ,
\end{eqnarray}
and by (\ref{a10}), (\ref{a14}), (\ref{m2}) and (\ref{f23}) for
all $\varepsilon, 0< \varepsilon <1,$ and $\eta >0$ we have
\begin{eqnarray}
\label{m10} &\log |f_\varepsilon ^3(e^{\eta +i\theta})|-\log
f_\varepsilon^3(e^\eta) )\leq \\ \nonumber & -C(e^{\eta/2}
+\sqrt{\varepsilon } e^{\eta /2}(\eta + \log \varepsilon )) + C
\leq -C b{'}(\eta, f_\varepsilon^3) +C, \quad \frac{\pi}{2}
<|\theta| < \pi .
\end{eqnarray}

Thus by (\ref{b28}) for all $f=f^1, f_\varepsilon^2, f_\varepsilon^3$
\begin{equation}\label{m11}
I_2(f) \leq C \exp ( - Ck), \quad k \geq k_0(f),
\end{equation}
where $k_0(f_\varepsilon^2),k_0( f_\varepsilon^3)$ do not depend
on $\varepsilon, 0 < \varepsilon <1.$

 From (\ref{b40}), (\ref{b41}), (\ref{m1}) and (\ref{m11})
we obtain Propositions 1 and 2.

{\bf Acknowledgments. } This work was written on the initiative of
\\ Prof.~ S.~Ju.~Favorov. The author sincerely thanks him for
support and help. The author is also deeply grateful to Dr.~
A.~M.~Vishnyakova for valuable discussions, important comments and
advice.

\end{document}